# A remark on $B_3$ − *sequences*

Li An-Ping

Beijing 100085, P.R. China
apli0001@sina.com

Abstract

A sequence of non-negative integers is called a $B_k$ − *sequences* if all the sums of arbitrary $k$ elements are different. In this paper, we will present a new upper bound for $B_3$ − *sequences*.

Keywords: $B_3$ − *sequences*, sums of arbitrary 3 elements, counting.

## 1. Introduction

Let $\mathcal{A} = \{a_i\}_1^r$, $0 \leq a_1 < a_2 < \cdots < a_r$, be a sequence of non-negative integers, $k$ a positive integer, if all the sums of arbitrary $k$ elements of $\mathcal{A}$ are different, we call $\mathcal{A}$ a $B_k - sequence$. Suppose $n$ is a arbitrary positive integer, denoted by $\Phi_k(n)$ the maximum of the sizes of the $B_k - sequences$ contained in $[0, n]$.

For $k = 2$, Erdos, Turan [3] and J. Singer [7] provided an upper bound and a lower bound respectively:

$$(1-\varepsilon)n^{1/2} \leq \Phi_2(n) \leq n^{1/2} + O(n^{1/4}). \tag{1.1}$$

B. Lindstrom [5] gave an improvement for the upper bound

$$\Phi_2(n) \leq n^{1/2} + n^{1/4} + 1. \tag{1.2}$$

For $k = 3$, S.C. Bose and S. Chowla [1] proved

$$\Phi_3(n) \geq (1-\varepsilon)n^{1/3},$$

and we presented an upper bound [4]

$$\Phi_3(n) \leq \left(\left(1 - \frac{1}{6\log_2^2 n}\right)4n\right)^{1/3} + 7. \tag{1.3}$$

For $k = 4$, B. Lindstrom [6] shown

$$\Phi_4(n) \leq (8n)^{1/4} + O(n^{1/8}), \tag{1.4}$$

J. Cilleruelo [2] gave a new upper bounds of $\Phi_3(n)$ and $\Phi_4(n)$

$$\Phi_3(n) \leq \left(\frac{4n}{1 + 16/(\pi+2)^4}\right)^{1/3} + o(n^{1/3}) \tag{1.5}$$

$$\Phi_4(n) \leq \left(\frac{8n}{1 + 16/(\pi+2)^4}\right)^{1/4} + o(n^{1/4}), \tag{1.6}$$

In this paper, we will show an improved upper bound over (1.3) for $\Phi_3(n)$

**Proposition 1.**

$$\Phi_3(n) \leq (3.962n)^{1/3} + 1.5. \tag{1.7}$$

## 2. The proof of main result

In this paper, it will be assumed $\mathfrak{S} = \{a_i\}_1^r$, $0 \leq a_1 < a_2 < \cdots < a_r \leq n$ is $B_3 - sequence$.

Denoted by $\lambda = \max\{(i-1.5)^3 / a_i \mid 1 \leq i \leq r\}$, and without loss the generality, we may assume that $\lambda \geq 3.962$, $a_r = n$, and $\lambda = (r-1.5)^3 / a_r$.

At first, we prove a simple lemma

**Lemma 1.** Suppose $a, b, c \in \mathfrak{S}$, $a \neq b \neq c$, and $x, y, z \in \mathfrak{S}$, $x \neq y \neq z$, if

$$a + b - c = x + y - z, \tag{2.1}$$

then $c = z$, $\{a, b\} = \{x, y\}$.

*Proof.* The equation (2.1) is rewritten as

$$a + b + z = x + y + c.$$

By the definition of $B_3 - sequence$, it follows

$$\{a, b, z\} = \{x, y, c\}.$$

And so,

$$c = z, \{a, b\} = \{x, y\}. \qquad \square$$

We introduce some symbols and notations will be used in this paper. Suppose that $X, Y, Z$ are three sets of integers in the interval $[0, n]$, define

$$\begin{aligned}
\mathfrak{I}_0(X, Y, Z) &= \{(x, y, z) \mid x \in X, y \in Y, z \in Z, x < y < z\}, \\
\mathfrak{I}_1(X, Y, Z) &= \{(x, y, z) \mid (x, y, z) \in \mathfrak{I}_0(X, Y, Z), x + y - z \geq 0\}, \\
\mathfrak{I}_2(X, Y, Z) &= \{(x, y, z) \mid (x, y, z) \in \mathfrak{I}_0(X, Y, Z), y + z - x < n\},
\end{aligned} \tag{2.2}$$

and

$$\begin{aligned}
\wp_0(X, Y, Z) &= \{x - y + z \mid (x, y, z) \in \mathfrak{I}(X, Y, Z)\}, \\
\wp_1(X, Y, Z) &= \{x + y - z \mid (x, y, z) \in \mathfrak{I}_1(X, Y, Z)\}, \\
\wp_2(X, Y, Z) &= \{y + z - x \mid (x, y, z) \in \mathfrak{I}_2(X, Y, Z)\}.
\end{aligned} \tag{2.3}$$

In the case $X = Y = Z$, we will simply write $\mathfrak{I}_i(X, Y, Z)$ and $\wp_i(X, Y, Z)$ as $\mathfrak{I}_i(X)$ and

$\wp_i(X)$ $i = 0,1,2$, respectively.

Besides, in the discussion we will take into account the other two classes of number sets

$$\wp_3(\mathfrak{S}) = \{2x - y \mid x, y \in \mathfrak{S}, 0 < 2x - y \leq n\},$$
$$\wp_4(\mathfrak{S}) = \{x \mid x \in \mathfrak{S}\}.$$

And define

$$\wp(\mathfrak{S}) = \bigcup_{0 \leq i \leq 4} \wp_i(\mathfrak{S})$$

By Lemma 1, we know that $|\wp_i(\mathfrak{S})| = |\mathfrak{I}_i(\mathfrak{S})|$, $i = 0,1,2$. and $\wp_i(\mathfrak{S}) \cap \wp_j(\mathfrak{S}) = \emptyset$, for $0 \leq i \neq j \leq 4$, hence,

$$|\wp(\mathfrak{S})| = \bigcup_{0 \leq i \leq 4} |\wp_i(\mathfrak{S})| \leq n. \tag{2.4}$$

Obviously,

$$|\wp_0(\mathfrak{S})| = \binom{r}{3}, \quad |\wp_4(\mathfrak{S})| = r. \tag{2.5}$$

Let $m = [n/2]$, $\mathfrak{S}_1 = \mathfrak{S} \cap [0, m], \mathfrak{S}_2 = \mathfrak{S} \cap (m, n]$, denoted by $|\mathfrak{S}_1| = s, |\mathfrak{S}_2| = t$.

For $a_i, a_j \in \mathfrak{S}_1, a_i < a_j$, then $0 < 2a_j - a_i \leq n$, and if $a_i, a_j \in \mathfrak{S}_2, a_i < a_j$, then $0 < 2a_i - a_j \leq n$, hence,

$$|\wp_3(\mathfrak{S})| \geq \binom{s}{2} + \binom{t}{2} \geq 2 \times \binom{(s+t)/2}{2} = \frac{r(r-2)}{4}. \tag{2.6}$$

In the next, we will mainly estimate $|\wp_1(\mathfrak{S})| + |\wp_2(\mathfrak{S})|$.

Denoted by $\mathcal{N} = \{0, 1, \ldots, n\}$, for $0 \leq i < j < k \leq m$, let

$$A(i, j, k) = \{(i, j, k), (i, j, n-k), (k, n-j, n-i), (n-k, n-j, n-i)\},$$
$$A'(i, j, k) = \{(i, k, n-j), (i, n-k, n-j), (j, n-k, n-i), (j, k, n-i)\},$$
$$B(i, j) = \{(i, i, j), (i, i, n-j), (j, n-i, n-i), (n-j, n-i, n-i)\},$$
$$B'(i, j) = \{(i, j, n-i), (i, n-j, n-i)\},$$
$$C(i, j) = \{(i, j, n-j), (j, n-j, n-i)\}.$$

It is easy to know that

$$(i, j, k), (i, j, n-k) \in \mathfrak{I}_2(\mathcal{N}), \quad (k, n-j, n-i), (n-k, n-j, n-i) \in \mathfrak{I}_1(\mathcal{N}).$$

That is,

$$A(i, j, k) \subseteq (\mathfrak{I}_1(\mathcal{N}) \cup \mathfrak{I}_2(\mathcal{N})).$$

Let
$$\mathcal{A} = \bigcup_{0 \le i < j < k \le m} A(i,j,k), \quad \mathcal{A}' = \bigcup_{0 \le i < j < k \le m} A'(i,j,k),$$
$$\mathcal{B} = \bigcup_{0 \le i < j \le m} B(i,j), \quad \mathcal{B}' = \bigcup_{0 \le i < j \le m} B'(i,j),$$
$$\mathcal{C} = \bigcup_{0 \le i < j \le m} C(i,j),$$

Suppose that $\chi(t)$ is the characteristic function of the set $\mathfrak{S}$, then it has
$$|\wp_i(\mathfrak{S})| = |\mathfrak{I}_i(\mathfrak{S})| = \sum_{(x,y,z) \in \mathfrak{I}_i(\mathcal{N})} \chi(x)\chi(y)\chi(z), \qquad i = 0,1,2.$$

**Lemma 2.** Let
$$\Gamma = \sum_{(x,y,z) \in \mathcal{A}} \chi(x)\chi(y)\chi(z),$$
$$\tilde{\Gamma} = \Gamma + \frac{1}{2} \sum_{(x,y,z) \in \mathcal{B}} \chi(x)\chi(y)\chi(z),$$
$$\overline{\Gamma} = \sum_{(x,y,z) \in \mathcal{A}' \cup \mathcal{B}'} \chi(x)\chi(y)\chi(z).$$

Denoted by $\Delta = \tilde{\Gamma} - \overline{\Gamma}$, then $\Delta \ge 0$, and
$$\Gamma \ge \frac{r^3}{12} + \frac{\Delta}{2} - \frac{3r^2}{8} - \frac{7r}{12}. \tag{2.7}$$

*Proof.*
$$\Delta = \frac{1}{2} \sum_{0 \le i < j < k \le m} (\chi(k) + \chi(n-k))(2\chi(i)\chi(j) + 2\chi(n-i)\chi(n-j))$$
$$+ \frac{1}{2} \sum_{0 \le i < k \le m} (\chi(k) + \chi(n-k))(\chi^2(i) + \chi^2(n-i))$$
$$- \frac{1}{2} \sum_{0 \le i < j < k \le m} (\chi(k) + \chi(n-k))(2\chi(i)\chi(n-j) + 2\chi(j)\chi(n-i))$$
$$= \sum_{k=0}^{m} \frac{\chi(k) + \chi(n-k)}{2} \left( \left( \sum_{0 \le v < k} \chi(v) \right)^2 + \left( \sum_{0 \le u < k} \chi(n-u) \right)^2 - 2 \sum_{0 \le v < k} \chi(v) \sum_{0 \le u < k} \chi(n-u) \right)$$
$$= \sum_{k=0}^{m} \frac{\chi(k) + \chi(n-k)}{2} \left( \sum_{0 \le v < k} (\chi(v) - \chi(n-v)) \right)^2 \ge 0.$$

On the other hand, it has
$$\frac{r^3}{6} \le \frac{1}{6} \times \left( \sum_{v=0}^{m} (\chi(v) + \chi(n-v)) \right)^3$$
$$\le \tilde{\Gamma} + \overline{\Gamma} + \frac{r^2}{4} + \frac{7r}{6}$$

$$\leq 2 \times \left(\Gamma + \frac{r^2}{4}\right) - \Delta + \frac{r^2}{4} + \frac{7r}{6}$$

$$\leq 2\Gamma - \Delta + \frac{3r^2}{4} + \frac{7r}{6}$$

i.e.

$$\Gamma \geq \frac{r^3}{12} + \frac{\Delta}{2} - \frac{3r^2}{8} - \frac{7r}{12}. \qquad \square$$

**Corollary 1.**

$$\Phi_3(n) \leq (4n)^{1/3} + 1. \tag{2.8}$$

*Proof.* It is easy to verify that (2.8) is true for $n \leq 6$, so we may assume that $r \geq 4$. Together with (2.4), (2.5), (2.6) and (2.7), it has

$$\binom{r}{3} + \frac{r^3}{12} + \frac{\Delta}{2} - \frac{3r^2}{8} - \frac{7r}{12} + \frac{r(r-2)}{4} + r - 1 \leq n$$

And

$$r^3 - 5r^2/2 + r \leq 4n.$$

So, for $r \geq 4$, it follows

$$(r-1)^3 \leq 4n. \qquad \square$$

Denoted by $\mathcal{A}[i]$ and $\mathcal{A}'[i]$ the $i$-th component of $\mathcal{A}$ and $\mathcal{A}'$ respectively, $i = 1, 2, 3, 4$, e.g. $\mathcal{A}[1] = \bigcup_{0 \leq i < j < k < m} (i, j, k)$. We know $(\mathcal{A}[1] \cup \mathcal{A}[2]) \subseteq \mathfrak{I}_2(\mathfrak{S})$, $(\mathcal{A}[3] \cup \mathcal{A}[4]) \subseteq \mathfrak{I}_1(\mathfrak{S})$.

Let

$$\mathcal{J} = (\mathfrak{I}_2(\mathfrak{S}) \setminus (\mathcal{A}[1] \cup \mathcal{A}[2])) \cup (\mathfrak{I}_1(\mathfrak{S}) \setminus (\mathcal{A}[3] \cup \mathcal{A}[4])),$$

Equivalency,

$$\mathcal{J} = (\mathfrak{I}_1(\mathfrak{S}) \cap (\mathcal{A}[1] \cup \mathcal{A}[2])) \cup (\mathfrak{I}_2(\mathfrak{S}) \cap (\mathcal{A}[3] \cup \mathcal{A}[4])) \cup ((\mathfrak{I}_1(\mathfrak{S}) \cup \mathfrak{I}_2(\mathfrak{S})) \cap \mathcal{A}').$$

We have the following estimation for $|\mathcal{J}|$

**Lemma 3.**

$$|\mathcal{J}| \geq \frac{(r-2.5)^3}{415.38} - \frac{r(r-1)}{4}. \tag{2.9}$$

The proof of Lemma 3 will be delayed in the later.

*The Proof of Proposition 1.* It is easy to verify directly that (1.7) is true for $n \leq 16$, so we may assume that $n \geq 17$, i.e. $r \geq 5$. Combine (2.4) ~ (2.9), it has

$$\binom{r}{3} + \frac{r^3}{12} + \frac{\Delta}{2} - \frac{3r^2}{8} - \frac{7r}{12} + \frac{r(r-2)}{4} + r - 1 + \frac{(r-2.5)^3}{415.38} - \frac{r(r-1)}{4} \leq n. \qquad (2.10)$$

It is easy to know that $\Delta \geq \frac{r-2}{2}$, hence

$$\frac{(r^3 - 3.54r^2 + 3.15r - 8.08)}{3.962} \leq n.$$

For $r \geq 5$, it follows

$$(r - 1.5)^3 \leq 3.962n.$$

The proof of Proposition 1 has been finished. □

In the next, we will prove (2.9).

For arbitrary integer $x, 0 \leq x \leq n$, define

$$\bar{x} = \begin{cases} x & \text{if } x \leq m \\ n - x & \text{if } x > m. \end{cases}$$

Define $\bar{\mathfrak{S}} = \{\bar{x} \mid x \in \mathfrak{S}\}$, $\bar{\mathfrak{S}}_2 = \{\bar{x} \mid x \in \mathfrak{S}_2\}$.

**Lemma 4.** For a positive integer $q, q \leq m$, denoted by $\mathfrak{A} = \mathfrak{S}_1 \cap [0, q], \mathfrak{B} = \bar{\mathfrak{S}}_2 \cap [0, q]$, $|\mathfrak{A}| = \alpha, |\mathfrak{B}| = \beta$. Let

$$\Lambda = |\wp_0(\mathfrak{A}, \mathfrak{A}, \mathfrak{B}) \cup \wp_0(\mathfrak{A}, \mathfrak{B}, \mathfrak{B}) \cup \wp_2(\mathfrak{A}, \mathfrak{A}, \mathfrak{B}) \cup \wp_1(\mathfrak{A}, \mathfrak{B}, \mathfrak{B})|,$$

then

$$\Lambda \geq \begin{cases} \dfrac{1}{4}\alpha\beta(\alpha + 2\beta - 4) - \dfrac{1}{12}\beta(\beta^2 - 4), & \text{if } \alpha \geq \beta, \\ \dfrac{1}{4}\alpha\beta(2\alpha + \beta - 4) - \dfrac{1}{12}\alpha(\alpha^2 - 4), & \text{if } \alpha < \beta. \end{cases} \qquad (2.11)$$

*Proof.* Suppose that $\mathfrak{A} = \{a_i\}_1^\alpha, a_1 < a_2 < \cdots < a_\alpha, \mathfrak{B} = \{b_i\}_1^\beta, b_1 < b_2 < \cdots < b_\beta$. For an integer $k, 1 \leq k \leq \beta$, denoted by

$$\mathfrak{A}_k = \{a_i \mid a_i + b_k < n, a_i \in \mathfrak{A}\}, \quad \mathfrak{A}'_k = \{a_i \mid a_i + b_k \geq n, a_i \in \mathfrak{A}\}.$$

Suppose $|\mathfrak{A}_k| = \alpha_k$, and $|\mathfrak{A}'_k| = \alpha - \alpha_k$, then there may be $\mathfrak{A}_k = \{a_i\}_1^{\alpha_k}$, $\mathfrak{A}'_k = \{a_i\}_{\alpha_k+1}^{\alpha}$. It is clear that for any two elements $a_i, a_j \in \mathfrak{A}_k, a_i < a_j$, then $(a_i, a_j, b_k) \in \mathfrak{I}_2(\mathfrak{A}, \mathfrak{A}, \mathfrak{B})$, and for any two elements $a_i, a_j \in \mathfrak{A}'_k, a_i < a_j$, then $(a_i, a_j, b_k) \in \mathfrak{I}_0(\mathfrak{A}, \mathfrak{A}, \mathfrak{B})$. For any one element $a_i \in \mathfrak{A}'_k$, and any two elements $b_u, b_v \in \mathfrak{B}, k \leq u < v$, then $(a_i, b_u, b_v) \in \mathfrak{I}_1(\mathfrak{A}, \mathfrak{B}, \mathfrak{B})$, and for any two elements $b_u, b_v \in \mathfrak{B}, 1 \leq u < v \leq k$, then $(a_i, b_u, b_v) \in \mathfrak{I}_0(\mathfrak{A}, \mathfrak{B}, \mathfrak{B})$. From the observation above, it follows

$$\Lambda \geq \sum_{i=1}^{\beta}\left(\binom{\alpha_i}{2} + \binom{\alpha - \alpha_i}{2}\right) + \sum_{i=0}^{\beta}(\alpha_i - \alpha_{i+1}) \cdot \left(\binom{\beta - i}{2} + \binom{i}{2}\right)$$

$$= \sum_{i=1}^{\beta}\left(\binom{\alpha_i}{2} + \binom{\alpha - \alpha_i}{2}\right) + \alpha \cdot \binom{\beta}{2} - \sum_{i=1}^{\beta}\alpha_i \cdot (\beta - 2i + 1).$$

If $\alpha \geq \beta$, then

$$\Lambda \geq \sum_{i=1}^{\alpha}\left(\binom{\left(\frac{\alpha+\beta}{2}\right) - i}{2} + \binom{\left(\frac{\alpha-\beta}{2}\right) + i}{2}\right) + \alpha \cdot \binom{\beta}{2} - \sum_{i=1}^{\beta}\left(\frac{\alpha+\beta}{2} - i\right) \cdot (\beta - 2i + 1)$$

$$\geq \alpha \cdot \binom{\beta}{2} + \frac{1}{4}\alpha\beta(\alpha - 2) - \frac{1}{12}\beta(\beta^2 - 4).$$

$$= \frac{1}{4}\alpha\beta(\alpha + 2\beta - 4) - \frac{1}{12}\beta(\beta^2 - 4).$$

else if $\beta \geq \alpha$, then

$$\Lambda \geq (\beta - \alpha) \cdot \binom{\alpha}{2} + \sum_{i=1}^{\alpha}\left(\binom{\alpha - i}{2} + \binom{i}{2}\right) - \sum_{j=1}^{(\beta-\alpha)/2}\alpha \cdot (\beta - 2j + 1) - \sum_{j=1}^{\alpha}(\alpha - j) \cdot (\alpha - 2j + 1)$$

$$\geq \beta \cdot \binom{\alpha}{2} + \frac{1}{4}\alpha\beta(\beta - 2) - \frac{1}{12}\alpha(\alpha - 1)(\alpha - 5).$$

$$\geq \frac{1}{4}\alpha\beta(2\alpha + \beta - 4) - \frac{1}{12}\alpha(\alpha^2 - 4).$$

□

**Corollary 2.** Suppose that $\overline{\mathfrak{S}} = \{z_i\}_1^r$, $z_1 \leq z_2 \leq \cdots \leq z_{r-1} \leq z_r \leq m$, then

$$z_k \geq \frac{7}{96}(k - 2)^3. \tag{2.12}$$

*Proof.* Let $q = z_k$ it is easy to know that $\left(\wp_0(\mathfrak{A}, \mathfrak{B}, \mathfrak{B}) \cup \wp_1(\mathfrak{A}, \mathfrak{B}, \mathfrak{B})\right) \subseteq [1, q]$, $\left(\wp_0(\mathfrak{A}, \mathfrak{A}, \mathfrak{B}) \cup \wp_2(\mathfrak{A}, \mathfrak{A}, \mathfrak{B})\right) \subseteq (n-q, n]$, and $\wp(\mathfrak{A}) \subseteq [1, q], \wp(\mathfrak{B}) \subseteq (n-q, n]$, so,

$$|\wp(\mathfrak{A})| + |\wp(\mathfrak{B})| + \Lambda \leq 2q.$$

By Corollary 1, we know

$$|\wp(\mathfrak{A})| \geq (\alpha-1)^3/4, |\wp(\mathfrak{B})| \geq (\beta-1)^3/4,$$

Hence, with Corollary 2,

$$\frac{1}{4}\alpha\beta(2\beta + \alpha - 4) - \frac{1}{12}\beta(\beta^2 - 4) + \frac{(\alpha-1)^3 + (\beta-1)^3}{4} \leq 2q, \quad \text{if } \alpha \geq \beta,$$

$$\frac{1}{4}\alpha\beta(2\alpha + \beta - 4) - \frac{1}{12}\alpha(\alpha^2 - 4) + \frac{(\alpha-1)^3 + (\beta-1)^3}{4} \leq 2q, \quad \text{if } \alpha < \beta.$$

With the symmetries, it follows

$$\frac{1}{4}\alpha\beta(\alpha + 2\beta - 4) - \frac{1}{12}\beta(\beta^2 - 4) + \frac{1}{4}\alpha\beta(2\alpha + \beta - 4) - \frac{1}{12}\alpha(\alpha^2 - 4) \leq 4q.$$

By the basic inequality, the left-hand side of the formula

$$\frac{1}{4}\alpha\beta(3\alpha + 3\beta - 8) + \frac{1}{12}(6(\alpha-1)^3 + 6(\beta-1)^3 - \alpha(\alpha^2-4) - \beta(\beta^2-4))$$

$$\geq \frac{7}{24}(\alpha + \beta - 2)^3.$$

That is,

$$z_k \geq \frac{7}{96}(\alpha + \beta - 2)^3 = \frac{7}{96}(k-2)^3.$$

$\square$

For a integers set $\mathcal{F}$ and a positive integer $q$, denoted by

$$\Delta_q(\mathcal{F}) = \left\{(x_i, x_j) \mid q < x_i, x_j \in \mathcal{F}, 0 < x_j - x_i \leq q,\right\}.$$

**Lemma 5.** Let $\overline{\mathfrak{S}}$ as defined before, $\overline{\mathfrak{S}} = \{z_i\}_1^r$, $z_1 \leq z_2 \leq \cdots \leq z_{r-1} \leq z_r \leq m$. Suppose that $q$ is a positive integer, $z_i \leq q < z_{i+1}$, then

$$\frac{(r-i)((q(r-i)/n)-1)}{2} \leq \left|\Delta_q(\mathfrak{S}_1) \cup \Delta_q(\overline{\mathfrak{S}}_2)\right|. \tag{2.13}$$

*Proof.* Let $d = [m/q]$, for each integer $t, t = 0, 1, \cdots, d$, denoted by $\mathfrak{D}_t = \overline{\mathfrak{S}} \cap [tq, (t+1)q)$,

$|\mathfrak{D}_t \cap \mathfrak{S}_1| = \mu_1(t)$, $|\mathfrak{D}_t \cap \bar{\mathfrak{S}}_2| = \mu_2(t)$. It is clear that for each $t, 1 \le t \le d$,

$$\binom{\mathfrak{D}_t \cap \mathfrak{S}_1}{2} \subseteq \Delta_q(\mathfrak{S}_1), \quad \binom{\mathfrak{D}_t \cap \bar{\mathfrak{S}}_2}{2} \subseteq \Delta_q(\bar{\mathfrak{S}}_2)$$

So,

$$\sum_t \left[ \binom{\mu_1(t)}{2} + \binom{\mu_2(t)}{2} \right] \le |\Delta_q(\mathfrak{S}_1) \cup \Delta_q(\bar{\mathfrak{S}}_2)|.$$

By the basic inequality,

$$\sum_t \left[ \binom{\mu_1(t)}{2} + \binom{\mu_2(t)}{2} \right] \ge \sum_t 2 \times \binom{(\mu_1(t) + \mu_2(t))/2}{2} = \sum_t 2 \times \binom{|\mathfrak{D}_t|/2}{2}$$

$$\ge 2d \times \binom{\sum_{1 \le t < d} |\mathfrak{D}_t| / 2d}{2}.$$

On the other hand,

$$\sum_{1 \le t < d} |\mathfrak{D}_t| = (r - i).$$

Hence,

$$|\Delta_q(\mathfrak{S}_1) \cup \Delta_q(\bar{\mathfrak{S}}_2)| \ge 2d \cdot \binom{(r-i)/2d}{2} = \frac{(r-i)(((r-i)/2d) - 1)}{2}$$

$$\ge \frac{(r-i)((q(r-i)/2m) - 1)}{2}$$

$\square$

*The Proof of Lemma 3.* For $q = z_i$, if $z_i \in \mathfrak{S}_1$, $(z_j, z_k) \in \Delta_{z_i}(\mathfrak{S}_1)$, then $(z_i, z_j, z_k) \in \mathfrak{I}_1(\mathfrak{S}_1)$; else if $(z_j, z_k) \in \Delta_{z_i}(\bar{\mathfrak{S}}_2)$, then $(z_i, n - z_k, n - z_j) \in \mathfrak{I}_1(\mathcal{A}'[2])$. If $z_i \in \bar{\mathfrak{S}}_2$, $(z_j, z_k) \in \Delta_{z_i}(\bar{\mathfrak{S}}_2)$, then $(n - z_k, n - z_j, n - z_i) \in \mathfrak{I}_2(\mathfrak{S}_2)$; else if $(z_j, z_k) \in \Delta_{z_i}(\mathfrak{S}_1)$, then $(z_j, z_k, n - z_i) \in \mathfrak{I}_2(\mathcal{A}'[4])$. Hence,

$$|\mathcal{J}| \ge \sum_{1 \le i \le r} |\Delta_{z_i}(\mathfrak{S}_1) \cup \Delta_{z_i}(\bar{\mathfrak{S}}_2)|.$$

Apply (2.12) and (2.13) by taking $q = z_i, i = 1, 2, \cdots, r,$ respectively, it follows

$$|\mathcal{J}| \ge \sum_{1 \le i \le r} \frac{(r-i) \times (z_i \cdot (r-i)/n - 1)}{2}$$

$$\ge \frac{7}{192n} \times \sum_{1 \le i \le r-2} i^3 \cdot (r-2-i)^2 - \frac{r(r-1)}{4}$$

$$\geq \frac{7(r-2)^2 \times ((r-2)^4 - 1)}{192 \times 60 \times n} - \frac{r(r-1)}{4}$$

$$\geq \frac{(r-2.5)^3}{415.38} - \frac{r(r-1)}{4}.$$

□